\documentclass[lettersize,journal]{article}

\usepackage{amsmath,amsfonts}
\usepackage{algorithmic}
\usepackage{algorithm}
\usepackage{xcolor}
\usepackage{array}
\usepackage[caption=false,font=normalsize,labelfont=sf,textfont=sf]{subfig}
\usepackage{textcomp}
\usepackage{stfloats}
\usepackage{hyperref}
\usepackage{url}
\usepackage{verbatim}
\usepackage{graphicx}
\usepackage{cite}
\usepackage{multirow} 
\usepackage{soul}
\usepackage{algorithm}
\usepackage{algorithmic}

\newcommand{\nonconverted}[1]{\mbox{}}
\usepackage{float}
\usepackage{amsthm}
\newtheorem{theorem}{Theorem}
\newtheorem{lemma}{Lemma}

\newtheorem{proposition}{Proposition}

\newcounter{broj}

\newcounter{broj1}
\def\dd{\mathrm{d}}

\newcommand{\te}{\textrm}
\newcommand{\tacka}{\,\cdot\,}

\newcommand{\veps}{\varepsilon}

\begin{document}

\title{Robust Frequency Domain Full-Waveform Inversion via HV-Geometry} 

\author{Zhijun Zeng, Matej Neumann, and Yunan Yang
\thanks{Zhijun Zeng is with the Department of Mathematical Sciences, Tsinghua University, Beijing 100084, China (e-mail: zengzj22@mails.tsinghua.edu.cn).}%
\thanks{Matej Neumann is with the Department of Mathematics, Cornell University, Ithaca, NY 14850, USA (e-mail: mn528@cornell.edu). Corresponding author.}%
\thanks{Yunan Yang is with the Department of Mathematics, Cornell University, Ithaca, NY 14850, USA (e-mail: yunan.yang@cornell.edu).}
\thanks{Zhijun Zeng and Matej Neumann contributed equally to this work.}}
\markboth{Journal of \LaTeX\ Class Files,~Vol.~14, No.~8, August~2021}%
{Shell \MakeLowercase{\textit{et al.}}: A Sample Article Using IEEEtran.cls for IEEE Journals}


\maketitle

\begin{abstract}
Conventional frequency-domain full-waveform inversion (FWI) is typically implemented with an $L^2$ misfit function, which suffers from challenges such as cycle skipping and sensitivity to noise. While the Wasserstein metric has proven effective in addressing these issues in time-domain FWI, its applicability in frequency-domain FWI is limited due to the complex-valued nature of the data and reduced transport-like dependency on wave speed. To mitigate these challenges, we introduce the HV metric ($d_{\text{HV}}$), inspired by optimal transport theory, which compares signals based on horizontal and vertical changes without requiring the normalization of data. We implement $d_{\text{HV}}$ as the misfit function in frequency-domain FWI and evaluate its performance on synthetic and real-world datasets from seismic imaging and ultrasound computed tomography (USCT). Numerical experiments demonstrate that $d_{\text{HV}}$ outperforms the $L^2$ and Wasserstein metrics in scenarios with limited prior model information and high noise while robustly improving inversion results on clinical USCT data.
\end{abstract}

\section{Introduction}
    

{Visualizing} and quantifying the acoustic properties of a medium from potentially sparse, indirect measurements constitutes the primary objective of imaging techniques across a broad spectrum of applications. Full-waveform inversion (FWI), originally introduced in~\cite{lailly1983seismic} and~\cite{tarantola1984inversion}, is a robust computational imaging method that has proven highly effective in both geophysical (seismic) and medical (ultrasound) applications. It operates by utilizing the entire recorded wavefield, incorporating both amplitude and phase information to iteratively refine a model of the medium. By minimizing the discrepancy between observed and simulated wave data, FWI facilitates the recovery of high-resolution characterizations of material properties, e.g., wave velocity or density.

FWI can be formulated in either the time-space domain \cite{tarantola1984inversion} or the frequency-space domain \cite{pratt1998gauss}; here we adopt the latter.  Frequency-domain FWI intrinsically allows the user to flexibly select discrete frequencies essential for a multiscale inversion strategy progressing hierarchically from low to high frequencies~\cite{pratt1999seismic}. This approach recovers large-scale structures from low frequencies and refines finer-scale details with higher frequencies, thereby mitigating the risk of being trapped in local minima. Additionally, frequency-domain FWI is sometimes computationally more efficient, as calculations are restricted to selected discrete frequencies, eliminating the need for time marching in the wave equation. Forward modeling in the frequency domain further allows optimal spatial discretization tailored to each frequency, with coarser grids for low frequencies and finer grids for high frequencies, maximizing accuracy while minimizing computational cost. Moreover, the frequency-domain formulation inherently accommodates attenuation by employing a complex-valued velocities, enabling accurate modeling of anelastic media\cite{malinowski2011high}.

In addition to its widespread application in seismic imaging, FWI has recently been adapted for medical imaging, particularly ultrasound computed tomography (USCT)\cite{lucka2021high,guasch2020full,ali2024ringFWI2D}. USCT utilizes an annular or cylindrical transducer array for data acquisition, sequentially emitting waves from each transducer and recording signals at the others, capturing both the transmitted and reflected waveforms from tissues. The acoustic field propagation is typically modeled using the wave equation in time-domain or the Fourier transform of the wave equation in frequency-domain, enabling FWI to achieve high-resolution reconstructions of tissue acoustic properties. Unlike seismic imaging, where FWI is constrained by partial boundary data at a fixed depth due to the half-space nature of the domain, USCT benefits from complete boundary data access across the spatial domain, including ideal transmission and refraction wave information.



Despite its widespread use, FWI continues to face significant computational challenges. The PDE-constrained optimization problem in FWI is inherently non-linear and non-convex, leading to being trapped in local minima. One such phenomenon is cycle-skipping, also known as phase wrapping in the frequency domain \cite{choi2015unwrapped}, which occurs when the synthetic data is misaligned with the observed data by more than half a wavelength, leading to incorrect velocity updates and trapping the inversion in local minima. The commonly used least-squares norm ($L^2$) as the objective function, traditionally applied in both time-domain and frequency-domain formulations,  is  sensitive to noise and susceptible to  cycle-skipping when paired with standard gradient descent methods. 
Extended FWI formulations based on the method of multipliers have demonstrated potential for mitigating cycle skipping by relaxing the wave‐equation constraint and introducing auxiliary wavefield variables that guide the inversion toward the global minimum~\cite{aghamiry2020robust,shah2014seismic,aghamiry2019improving,gholami2024full,ranjbaran2023quantitative}. However, these methods entail a substantial increase in computational cost and algorithmic complexity.  Another way to mitigate cycle skipping is to design alternative misfit functions; in recent years, numerous studies have proposed such functions, based on instantaneous phase, adaptive matching filters, and envelope measurements, to improve the inversion objective's convexity and enhance robustness to noise~\cite{bozdaug2011misfit,wu2014seismic,zhu2016building}. More recently, the Wasserstein metric, originating from optimal transport theory~\cite{villani2021topics}, has emerged as a robust framework for comparing seismic signals \cite{engquist2013application,engquist2016optimal,yang2017analysis,yang2018application,engquist2020quadratic,engquist2022optimal} and tackling cycle-skipping issues. However, foundational differences in optimal transport theory, specifically the focus on probability distributions, contrast with the properties of seismic signals, which are often signed or complex-valued. As a result, the full potential of optimal transport-based techniques remains largely untapped for frequency-domain FWI.

\begin{figure*}[htbp!]
\centering
\includegraphics[width=\textwidth]{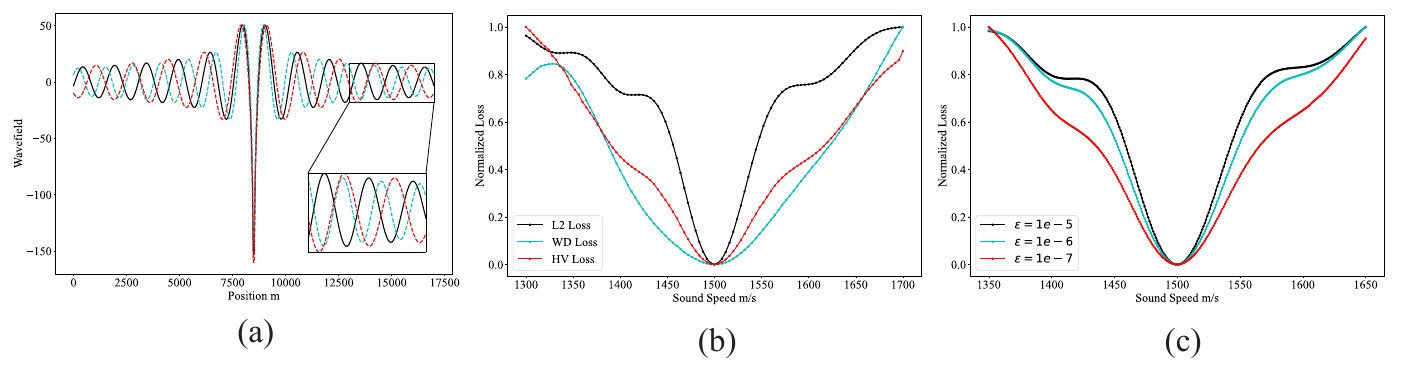}%
\caption{(a) Simulated frequency-domain data is generated for the wave speed $c(x)$, ranging from 1300 to 1700$\,\mathrm{m/s}$. (b) Misfit functions are computed by comparing reference data, based on the true wave speed $c^* = 1500 \,\mathrm{m/s} $, with simulated data corresponding to different wave speeds $ c(x) = c \in (1300, 1700) \,\mathrm{m/s} $. Three misfit functions are evaluated: the squared $L^2$ norm, the squared 2-Wasserstein metric, and the squared HV metric. (c) {HV metric between the reference data, based on the true wave speed $c^* = 1500 \,\mathrm{m/s} $, and the simulated data corresponding to different wave speeds $ c(x) = c \in (1350, 1650) \,\mathrm{m/s} $, are computed. Three HV metric with different $\varepsilon\in\{10^{-5}, 10^{-6}, 10^{-7}\}$ are evaluated.}}
\label{fig_sim}
\end{figure*}

Similarly to the  $L^2$ misfit,  the Wasserstein distance  can also face cycle-skipping issues in both frequency-domain and time domain FWI  due to data normalization ~\cite[Fig.~2]{yang2018application}. For an illustration, we adopt a constant‑velocity model with $c(x)=c$, where $c$ varies from 1300 to 1700 m/s. Twenty sources and 851 receivers are uniformly spaced at a depth of 100m, and a 1Hz signal is employed (see Section~\ref{sec:noise_free_seismic} for details).
Fig.~\ref{fig_sim}(a) illustrates the observed frequency-domain signals for various $c$ values in seismic imaging, revealing phase shifts in signals recorded at locations away from the source. Furthermore, in Fig.~\ref{fig_sim}(b), we compute the $L^2$ misfit and the Wasserstein metric (normalized via linear transformation~\cite{Survey2,engquist2022optimal}) between the signals corresponding to $ c(x)=c $ over the range $[1300, 1700]\,\mathrm{m/s}$ and those corresponding to the $ c(x)=1500\,\mathrm{m/s} $ model. The squared $L^2$ norm exhibits pronounced nonconvexity. The squared Wasserstein metric has a much wider basin of attraction than the $L^2$ norm. However, it becomes non-monotonic for shifts greater than $200\,\mathrm{m/s}$, which illustrates the potential limitation of the Wasserstein metric in handling convexity issues of frequency-domain FWI. These observations motivate us to consider a different objective function for frequency-domain FWI.

To characterize rich features of real-world signals, Miller et al.~\cite{miller2001group} introduced a novel geometric framework that accommodates sign-changing signals and quantifies both horizontal deformations (measured by the Wasserstein metric) and vertical deformations (characterized by the $L^2$ norm). In~\cite{han2023hv}, further developments were made to this framework, including characterizing the regularity and stability of minimizing geodesics and introducing an efficient fixed-point iteration-based algorithm validated by several numerical examples. Given a finite spatial interval (assumed to be $[0,1]$ for simplicity), a path $f(x,t)$ connecting two one-dimensional signals, $f_0$ and $f_1$, is described by the solutions of the transport equation with time-dependent source term $z(x,t)$ and velocity $v(x,t)$:
\begin{equation}\label{eq:adm}
\begin{aligned}
f_t &= -f_x v + z, \quad \text{on } [0,1] \times [0,1], \\
v(0,\cdot) &= v(1,\cdot) = 0, \\
f(\cdot,0) &= f_0, \quad f(\cdot,1) = f_1.
\end{aligned}
\end{equation}
We denote the set of all tuples $(f,v,z)$ satisfying~\eqref{eq:adm} as the admissible set $\mathcal{A}$. For a fixed path represented by $(f,v,z)$ and hyperparameters $(\kappa, \lambda, \veps)$, the action measuring the effort of deforming $f_0$ to $f_1$ along the particular path is 
\begin{equation} \label{eq:action1}
   A_{\kappa, \lambda, \veps}(f,v,z) =\frac12 \int_0^1 \sqrt{\int_0^1   \kappa v^2 + \lambda v_x^2  + \varepsilon v_{xx}^2 + z^2    \, \dd x} \dd t
\end{equation}
where $\kappa >0$, $\lambda \geq 0$, and $\veps>0$. The cost of the action highly depends on the choice of hyperparameters $(\kappa,\lambda,\varepsilon)$ as seen in Fig.~\ref{fig:convexpara}. Also see examples in \cite{han2023hv}.
The $v^2$ term measures the horizontal movement (modeled by the transport equation term $\partial_x f \cdot v$) of the signal, while the $z^2$ term measures the vertical change. The hyperparameter $\kappa$ is present for modeling flexibility, as it accounts for horizontal and vertical variations differently, while the $\veps v_{xx}^2$ term is necessary for regularity (see details in~\cite{han2023hv}). 
Using the action we can define the HV metric between functions $f_0$ and $f_1$ as the infimum of the action over the set $\mathcal{A}$:
\begin{equation}\label{eq:dhv}
\begin{aligned}
d_{HV(\kappa, \lambda, \varepsilon)} &:= 
 \inf_{(f,v,z)\in \mathcal{A}} \sqrt{A_{\kappa, \lambda, \varepsilon}(f,v,z)}.
\end{aligned}
\end{equation}

The HV metric can also be viewed as an extension of optimal transport metrics. Table~\ref{tab:comparison} provides a comparison among the HV metric~\cite{han2023hv}, the dynamic formulation of the 2-Wasserstein ($W_2$) metric from optimal transport (OT)~\cite{villani2021topics}, and unbalanced optimal transport (UOT)~\cite{chizat2018scaling}. All these metrics can be interpreted as the optimal cost computed over a set of admissible paths; however, their distinctions lie in the types of paths considered and the specific cost functionals used.

\begin{table}[!ht]
\begin{tabular}{c|l}
\hline
\multicolumn{1}{c|}{Metric}                            & Formulation \\ \hline
HV& $\left\{\begin{array}{lll}
\partial_t f &=& -f_x v + z\\
A_{\text{HV}}&=&\frac12 \int_0^1 \sqrt{\int_0^1 \left(\kappa v^2 + \lambda v_x^2 + \varepsilon v_{xx}^2 + z^2 \right)  \dd x} \dd t 
\end{array}\right.$
\\ \hline
OT& $\left\{ \begin{array}{lll}
\partial_t f &=& -\text{div}(fv)\\
A_{\text{OT}} &=& \int_0^1\int_0^1 v^2 f \dd x\dd t 
\end{array}\right.$
\\ \hline
UOT& $\left\{ \begin{array}{lll}
\partial_t f &=& -\text{div}(fv)+zf\\
A_{\text{UOT}} &=& \int_0^1\int_0^1 (v^2+z^2) f \dd x\dd t
\end{array}\right.$
\\ \hline
\end{tabular}
\caption{Comparison of the HV metric, optimal transport metric with quadratic cost, and unbalanced optimal transport metric.}
\label{tab:comparison}
\end{table}

In this study, we propose a novel misfit function---the HV metric~\cite{han2023hv}---for frequency-domain FWI. Unlike the Wasserstein metric, the HV metric naturally extends to complex-valued signals while retaining its transport-like properties. Importantly, this extension supports signed signals without requiring equal total mass. It eliminates the need for pre-processing or artificial interventions to meet mathematical prerequisites.  

For the numerical implementation of $d_{\text{HV}}$, \cite{han2023hv} proposed an alternating direction optimization algorithm to solve the variational problem \eqref{eq:dhv} for one-dimensional real-valued signals. In the case of complex-valued signals in the frequency domain, we adopt a product measure to facilitate a straightforward calculation. Fig.~\ref{fig_sim}(b) illustrates the corresponding loss curve for a single-parameter model, where the HV metric demonstrates smooth and monotonic behavior, better than both the $L^2$ and the Wasserstein ones. Notably, it avoids the staircase patterns commonly observed in other metrics, which reflects its potential as a robust and effective misfit function for frequency-domain FWI.


The remainder of this paper is organized as follows. In Section~\ref{sec:HV}, we introduce the definition and key properties of $d_{\text{HV}}$. We then present a finite-difference-based alternating direction optimization method for numerically computing $d_{\text{HV}}$, along with the calculation of its Fréchet derivative for complex-valued signals. Section~\ref{sec:ne} showcases three numerical experiments: first, we evaluate the performance of the proposed HV metric in seismic inversion using the Marmousi model, a classic seismic benchmark, under both noise-free and noisy conditions; second, we compare the effectiveness of various metrics, as illustrated in Fig.~\ref{fig_sim}, for USCT imaging of a synthetic breast model with single-band noisy data; third, we demonstrate the robustness and effectiveness of the HV metric on real-world data for USCT imaging using multi-band clinical breast data. These studies reveal that the HV metric outperforms the conventional $L^2$ norm and Wasserstein metric in mitigating cycle-skipping issues and exhibits superior noise robustness in frequency-domain FWI. Finally, Section~\ref{sec:conclusion} summarizes our findings and concludes the study.

\section{Background}\label{sec:HV}
In this section, we introduce some essential mathematical background on the HV metric.
\subsection{Frequency‐Domain FWI via Adjoint‐State Method}
 Mathematically, the FWI problem can be formulated as a nonlinear Partial Differential Equation (PDE)-constrained optimization problem. In the frequency domain, it seeks to, for a fixed frequency $\omega,$ find the velocity model $c(\mathbf{x})$ by minimizing a general misfit
\begin{eqnarray}
  & \min_{c}   &  J(c)
  =
  \sum_{k=1}^K \mathcal{D} \bigl(u_k(\mathbf{x}_r,\omega),\,d_k^{\rm obs}(\mathbf{x}_r,\omega;c)\bigr) \label{eq:Misfit} \\
  &\text{s.t.}   & \Bigl[\nabla^2 + \tfrac{\omega^2}{c^2(\mathbf{x})}\Bigr]\,u_k(\mathbf{x})
  = -\,s_k(\mathbf{x}), \label{eq:helm}
\end{eqnarray}
where $K$ is the total number of sources, $s_k$ is the $k$-th source term, $u_k$ is the synthetic wavefield with the source $s_k$, $d_k^{\rm obs}$ is the observed data with the source $s_k$, $\mathcal{D}$ is any metric comparing two datasets (e.g., $L^2$ and the HV metric), and $\mathbf{x}_r$ is the transducer locations. Equation~\eqref{eq:helm} is the Helmholtz equation which implicitly decides the forward operator in the inverse problem: $c\mapsto \{u_k\}_{k=1}^K$.

The adjoint method is a widely-used approach for numerical solution of the PDE-constrained optimization problem \eqref{eq:Misfit}-\eqref{eq:helm}. Introducing Lagrange multipliers $\{\lambda_k(\mathbf{x})\}_{k=1}^K$, the Lagrangian is formed as
\begin{equation}
  \mathcal{L}
  = J(c)
-
  \sum_{k=1}^K
  \Re\bigl\langle
    \lambda_k,\,
    \bigl[\nabla^2 + \tfrac{\omega^2}{c^2}\bigr]\,u_k + s_k
  \bigr\rangle,
  \label{eq:Lagrangian}
\end{equation}
with $\langle f,g\rangle=\int_\Omega f^*\,g\,\mathrm{d}\mathbf{x}$. By requiring the first-order optimality $\delta\mathcal{L}/\delta u_k=0$ we obtain the adjoint equation for $\lambda_k$:
\begin{equation}
  \Bigl[\nabla^2 + \tfrac{\omega^2}{c^2(\mathbf{x})}\Bigr]\lambda_k(\mathbf{x})
  =   \frac{\delta \mathcal{D}}{\delta u_k}\left(u_k(\mathbf{x}_r),\,d_k^{\rm obs}(\mathbf{x}_r)\right)
  \delta(\mathbf{x}-\mathbf{x}_r),
  \label{eq:Adjoint}
\end{equation}
where the right‐hand side is the variational derivative of the misfit with respect to $u_k$, injected at receiver positions. The gradient of $J$ depends on the solutions $\{\lambda_k\}_{k=1}^K$ to the adjoint equations~\eqref{eq:Adjoint}:
\begin{align}
  \frac{\delta J}{\delta c}(\mathbf{x}) \nonumber
  =& -\sum_{k=1}^K 
    \Re\Bigl\langle
      \lambda_k,\,
      \frac{\partial}{\partial c}\!\Bigl[\nabla^2 + \tfrac{\omega^2}{c^2}\Bigr]\,u_k
    \Bigr\rangle\\
  =& -\,2\omega^2\sum_{k=1}^K
    \Re\frac{\lambda_k^*(\mathbf{x})\,u_k(\mathbf{x})}{c^3(\mathbf{x})}.
  \label{eq:Gradient}
\end{align}
We address that a different choice of the metric $\mathcal{D}$ only affects the source term in the adjoint equation~\eqref{eq:Adjoint}. To obtain the gradient, we need to solve $K$ forward wave propagation~\eqref{eq:helm} and $K$ adjoint equations~\eqref{eq:Adjoint} once the adjoint sources are computed.
\begin{algorithm}[htbp]
  \caption{Frequency‐Domain FWI via Adjoint Method}
  \label{alg:fd-fwi}
  \begin{algorithmic}[1]
    \REQUIRE initial model $c^{(0)}$, data $\{d_k^{\rm obs}\}$, frequencies $\{\omega\}$
    \FOR{$n = 0,1,2,\dots$ \textbf{until} convergence}
      \FOR{each source $k$ and frequency $\omega$}
        \STATE Solve forward: 
        $
          \bigl[\nabla^2 + \tfrac{\omega^2}{(c^{(n)})^2}\bigr]\,u_k^{(n)} = -\,s_k
        $
        \STATE Compute residual:
        $
          r_k^{(n)}
          = \frac{\partial\mathcal{D}}{\partial u_k}\bigl(u_k^{(n)},d_k^{\rm obs}\bigr)
        $
        \STATE Solve adjoint:
        $
          \bigl[\nabla^2 + \tfrac{\omega^2}{(c^{(n)})^2}\bigr]\,\lambda_k^{(n)}
          = r_k^{(n)}\,\delta(\mathbf{x}-\mathbf{x}_r)
        $
      \ENDFOR
      \STATE Assemble gradient $g^{(n)}(\mathbf{x})$ via \eqref{eq:Gradient}.
      \STATE Update model:
      $
        c^{(n+1)} = \text{Optim}(c^{(n)},,g^{(n)})
      $
    \ENDFOR
  \end{algorithmic}
\end{algorithm}
 
\subsection{The HV metric and its properties}
For two signal function $f_0,f_1 \in L^2(0,1)$ and the space  $\mathcal V:= L^2((0,1),H^2(0,1) \cap H^1_0(0,1))$, we define the set of all admissible paths $\mathcal A(f_0, f_1)$ as all paths satisfying \eqref{eq:adm} weakly (see \cite{han2023hv} for rigorous definition)
and the HV metric by Eqn.~\eqref{eq:dhv}. Under the aforementioned definitions, \cite{han2023hv} conducted a detailed analysis of $d_{\text{HV}}$, demonstrating~\eqref{eq:dhv} indeed defines a metric together with several properties described in Theorem~\ref{thm:prop} below.

\begin{theorem}\label{thm:prop} The action and $d_{\text{HV}}$ have the following properties:
\begin{enumerate}
    \item Consider  $f_0,f_1\in L^{2}(0,1)$. 
There exists an admissible path $(f,v,z) \in \mathcal A(f_0, f_1) $ minimizing the action \eqref{eq:action1}.
\item $d_{\text{HV}}$ is a metric on $L^{2}(0,1)$.
\item The metric space $(L^2(0,1),d_{\text{HV}})$ is complete.
\end{enumerate}    
\end{theorem}
Given $f_0,f_1 \in H^1(0,1)$, solving the minimizing problem in $d_{\text{HV}}$
requires deriving the corresponding first-order optimality conditions, commonly called the Euler--Lagrange equations. These conditions are described as follows:
\begin{proposition}
    If $f_0, f_1 \in H^1(0,1),$ then the minimizing path $f \in L^2((0,1),H^1(0,1))$, and the resulting 
  Euler--Lagrange equations satisfied by $f,v,z$, combined with the condition for belonging to $\mathcal A(f_0, f_1)$ can be expressed as follows:
\begin{align}
\veps v_{xxxx} -\lambda v_{xx} + \kappa v + z f_x & = 0 \te{ weakly on } (0,1)^2 \label{eq:elv}\\
 v = 0   \; \te{ and } \;  v_{xx} & = 0 \,\,\te{ on } \partial \Omega \times [0,1] \;  \label{eq:v-BCs}\\
z_t + (zv)_x & = 0  \te{ weakly on } (0,1)^2 \label{eq:elz}  \\
 f_t + f_x  v - z & = 0 \te{ weakly on } (0,1)^2 \label{eq:elf}\\
 f(\tacka,0)  =f_0, \; \te{ and }  \; f(\tacka,1)&=f_1. \label{eq:f-ICs}
\end{align}
\end{proposition}

\subsection{Computation of the HV metric}
This coupled PDE system above is complex and challenging to solve directly. To obtain its numerical solution, \cite{han2023hv} splits the system of equations into two sub-problems, each of which can be regarded as a convex optimization problem, and performs an alternating-direction optimization algorithm to iteratively solve both problems. 
We will also provide both the flops and memory complexity for each step of the algorithm. \medskip

\subsubsection{Fixed $v$, optimize $(f,z)$}
For any given $v$, the optimization on $(f,z)$ in the definition of the HV metric reduces to a quadratic objective functional under a linear constraint
\begin{equation}\label{eq:sub2}
  \min_{f,z} \frac12 \int_0^1 \int_0^1  z^2  \, \dd x \dd t,  \quad \text{s.t.}\quad (f,v,z) \in \mathcal{A},
\end{equation}
where the first-order optimality conditions are given by equations:
\begin{align}\label{eq:sub1}
\begin{split}
z_t + (zv)_x & = 0,  \\
f_t + f_x  v - z & = 0,  \\
f(\tacka,0) =f_0,\,
\,& f(\tacka,1)=f_1. 
\end{split}
\end{align}
In \cite{han2023hv}, the authors adopt a Lagrange approach to solve~\eqref{eq:sub1} with the following mathematical guarantee. 
\begin{lemma}\label{lem:representation_sol}
 Given $v\in\mathcal V$, $z\in L^2((0,1)^2)$, let $f_0\in L^2(0,1)$. If $f\in L^{2}((0,1)^2)$ is a weak solution to the initial value problem
 \[
 f(\cdot,0)=f_0,\quad\partial_t f+\partial_x f\cdot v=z \quad \;\te{on } [0,1]^2,
 \]
 then
$f$ has the following representation: for almost every $x,t$ in $[0,1]^2$,
\begin{equation}\label{eq:solf}
f(\Phi(x,t),t)=f_0(x)+\int_0^t z(\Phi(x,s),s)\,\dd s, 
\end{equation}
where $\Phi$ is the flow of the vector field $v$:
\begin{equation}\label{eq:flow-map}
\partial_t \Phi(x,t)=v(\Phi(t,x),t),\quad \Phi(0,x)=x.
\end{equation}
It follows that the weak solution $f$ is unique. 
 \end{lemma}
As a result, the solution to~Eqn.~\eqref{eq:elz}  can be expressed as
\begin{equation} \label{eq:solz}
z( \Phi ( x, t), t) = z( x, 0 )\, e^{-\int_0^t v_x ( \Phi(x, s), s) \,\mathrm{d}s},
\end{equation}
where $\Phi$ denotes the flow of $v$, appearing in Lemma~\ref{lem:representation_sol}. Introducing an auxiliary variable $J(x,t) = e^{-\int_0^t v_x ( \Phi(x, s), s) \,\mathrm{d}s}$, we obtain the representation of $f$ as:
\[
\begin{aligned}
f_1(\Phi(x,1)) &= f_0(x) + \int_0^1 z( \Phi(x, t), t)\,\mathrm{d}t\\ &= f_0(x) + z(x,0) \int_0^1 J(x,t) \,\mathrm{d}t.
\end{aligned}
\]
Using the boundary conditions, we finally derive the following representation formulas for action-minimizing paths:
\begin{proposition}\label{prop:v2zf}
    Given $v\in\mathcal V$, the weak solution $(z,f)$ of \eqref{eq:f-ICs} has the following analytical representation
    \begin{equation*}
     \begin{aligned}
        z( \Phi ( x, t) , t) &= \left( f_1 ( \Phi(x, 1) ) - f_0 (x) \right) \, \frac{J(x,t)}{\int_{0}^{1} J(x,\tau)  \dd \tau}, \\
        f(\Phi(x,t),t) &= \left(1-\eta(x,t) \right) \,  f_0( x) +  \eta(x,t) \, f_1 ( \Phi ( x, 1)), \\\te{where } \eta(x,t) &=  \textstyle\int_0^t J(x,s) \dd s\, /\textstyle\int_0^1 J(x,s) \dd s.
    \end{aligned}
    \end{equation*}
\end{proposition}

By solving Eqn.~\eqref{eq:flow-map}, we obtain $\Phi(x,t)$ in $O(N_xN_t)$ flops and memory where $N_x$ and $N_t$ represent the number of spatial and time grid points, respectively. The analytical solutions for $z(\Phi(x,t), t)$ and $f(\Phi(x,t), t)$ provided in  Prop.~\ref{prop:v2zf}  are computed using first-order numerical integration. Using these results, we then obtain the functions $f$ and $z$ at the desired $(x,t)$ coordinates using linear interpolation. \medskip

\subsubsection{Fixed $f$, optimize $(v, z)$}
\label{timecomplex}
Given $f$ from the previous step, the second subproblem involves optimizing the pair $(v, z)$ to minimize the objective function while satisfying the same constraints. The Euler–Lagrange equation, combined with the constraint $z = f_t + v f_x$, yields a fourth-order boundary value problem for $v$:
\begin{align}\label{eq:sub3}
\begin{split}
\varepsilon v_{xxxx} -\lambda v_{xx} + (\kappa  + |f_x|^2 ) v  & = -f_t \, f_x  \quad \text{on } (0,1)^2, \\
v & = 0  \quad \text{at } x = 0,1,\ \forall t, \\
v_{xx} & = 0  \quad \text{at } x = 0,1,\ \forall t.
\end{split}
\end{align}
Once obtaining $v$, the function $z$ is computed as $z = f_t + v f_x$.

For the system \eqref{eq:sub3}, fixing $ t $ reduces the problem to a fourth-order ordinary differential equation (ODE) for $v(\cdot, t)$ with boundary conditions $v = 0$ and $v_{xx} = 0$. This can be solved using the finite-difference method on a uniform mesh over the space-time domain $[0,1] \times [0,1]$. Letting $\Delta x = \frac{1}{N_x}$ and $\Delta t = \frac{1}{N_t}$ denote the spatial and temporal grid spacings, respectively, we can express the discrete form of \eqref{eq:sub3} as:
\begin{equation}\label{eq:lin_equation}
\left(A+\operatorname{diag}\left(\mathbf{w}_{j} \odot \mathbf{w}_{j}\right)\right) \mathbf{v}_{j}=-\boldsymbol{\tau}_{j} \odot \mathbf{w}_{j},
\end{equation}
where $\operatorname{diag}(\cdot)$ represents a diagonal matrix with elements of the vector $\cdot$ as its diagonal entries, $A \in \mathbb{R}^{(N_x+1) \times (N_x+1)}$, and $\mathbf{w}_{j}$ and $\boldsymbol{\tau}_{j}$ as defined in \cite{han2023hv}.  Similarly to the complexity of the previous subsection, the memory and flops complexity of creating the matrix is $O(N_xN_t)$. Since the matrix $A$ in (\ref{eq:lin_equation}) is sparse with bandwidth equal to 5, solving for $v_j$ requires $O(N_x)$ memory and flops. The flops complexity of computing the full HV metric is equal to $O(N_xN_tN_i)$, where $N_x$ and $N_t$ represent the number of space and time grid points, and $N_i$ represents the number of iterations we perform the algorithm. The memory complexity is equal to $O(N_xN_t)$ since at each point we only save the current values of $f,v$, and $z$. We remark that in all applications, $N_t$ is usually small compared to $N_x$, simplifying the flops complexity to $O(N_xN_i)$ and memory complexity to $O(N_x)$.  \medskip

\subsection{Dependence on hyperparameters $\kappa,\lambda$ and $\varepsilon$} 
\begin{figure}[!ht]
    \centering    \includegraphics[width=0.5\linewidth]{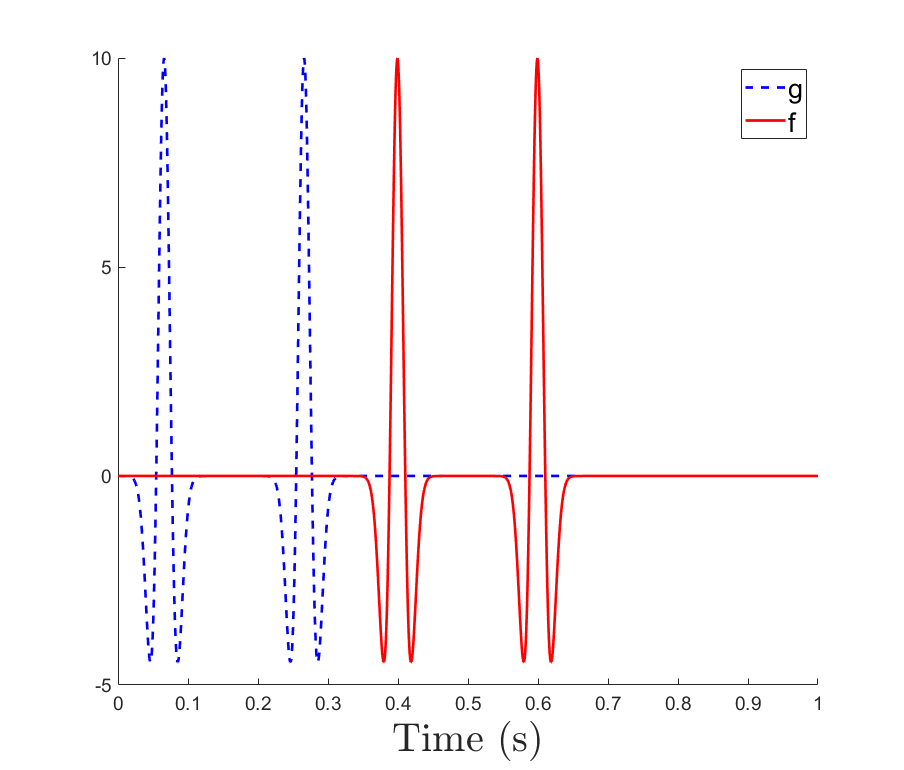}
    \caption{Signal $f$ and the shifted signal $g\left(t;\frac{1}{3}\right)=f\left(t-\frac13\right)$.}
    \label{Ricker}
\end{figure}

Eqn.~\eqref{eq:action1} suggests a certain dependence of the HV metric on the choice of hyperparameters $\kappa,\lambda$ and $\varepsilon$. By increasing $\kappa$, we penalize velocity, making vertical deformations less costly and making the HV metric behave similarly to $L^2$. On the other hand, decreasing the hyperparameters lowers the cost of horizontal deformations, thus encouraging them, and making the HV metric more transport-like. This apparent interpolation feature between vertical and horizontal deformations makes the HV metric useful in many applications, as it can be adapted to a wide range of problems. We showcase this property by focusing on two pairs of two Ricker wavelet signals, where one pair is a time shift of the other. We regard them as observed data $f(t)$ and synthetic data $g(t;s)=f(t-s)$, respectively, with $s$ quantifying the shift in time as shown in Fig.~\ref{Ricker}.

\begin{figure}[!ht]
    \centering
    \includegraphics[width=1\linewidth]{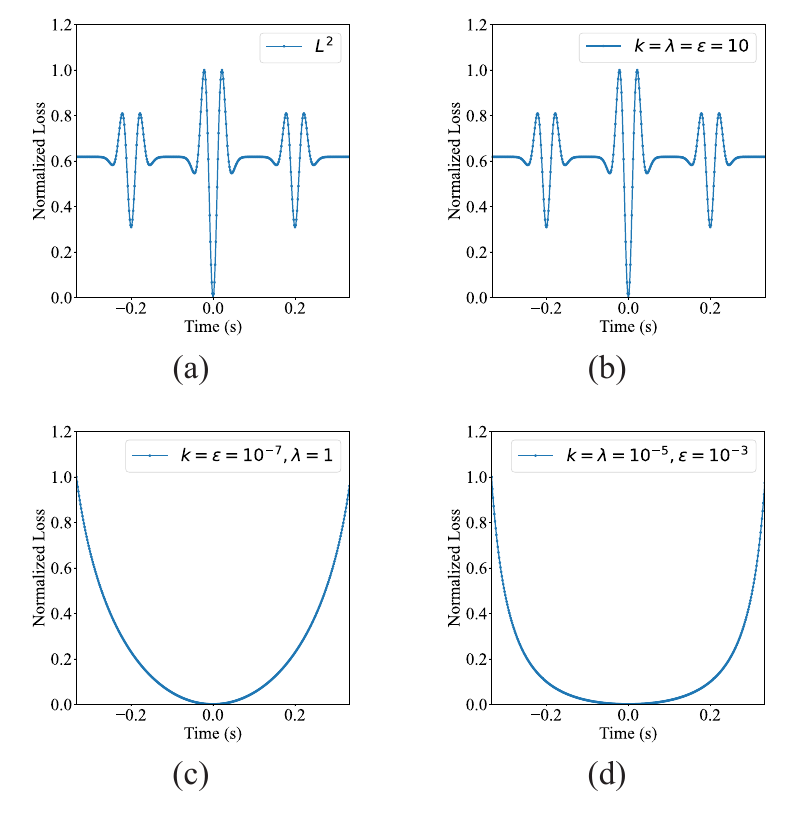}
    \caption{The normalized loss between $f(t)$ and $g(t) = f(t-s)$ using the: (a) $L^2$ norm; (b) HV metric with $\kappa=\lambda=\varepsilon=10$; (c) HV metric with $\kappa=\varepsilon=10^{-7}$ and $\lambda=1$; (d) HV metric with $\kappa=\lambda=10^{-5}$ and $\varepsilon=10^{-3}$.}
    \label{fig:convexpara}
\end{figure}

We plot the $L^2$ norm between $f$ and $g$ alongside the HV metric results with different choices of hyperparameters in Fig.~\ref{fig:convexpara}. Also see Fig.~\ref{fig_sim} for a similar frequency domain analysis. We observe that both HV and $L^2$ behave similarly when $\kappa$ is large and displays many local minima and maxima, demonstrating the difficulty of the cycle-skipping issues in time-domain FWI. As we decrease the hyperparameters, the local minima and maxima disappear, sharing similar results with the 2-Wasserstein metric; see \cite{Survey2}. Decreasing the hyperparameters even further, we notice that the plot becomes flatter around the global minimum, suggesting that the HV metric, with very small hyperparameters, shares features with weaker norms than the Wasserstein metric and its equivalent $\dot{H}^{-1}$ norm~\cite{otto2000generalization}.

\subsection{Fr\'echet derivative of the HV metric}
Given two signals $f_0,f_1\in L^2(0,1)$ let the tuple $(f,v,z)$ denote the optimal admissible path between $f_0$ and $f_1$. Next, we compute the Fr\'echet derivative of the squared HV metric with respect to $f_0$, i.e., the functional 
$$
J: f_0\rightarrow d_{\text{HV}}^2(f_0, f_1).
$$
This is an important quantity in any gradient-based optimization algorithms when $d_{\text{HV}}^2$ is used as the data-matching misfit function, including steepest descent, stochastic gradient descent, and quasi-Newton methods. 

To proceed, we follow the typical steps to obtain the directional derivative of the functional. Consider a curve $g_s(x)$ in $L^2[0,1]$ parameterized by $s$ with $f_0(x) =  g_{s=0}(x)$. Let $(w,h)$ denote the HV geometry tangent vectors in the tangent space at point $f_0$ corresponding to the Eulerian perturbation $\partial_s g_s|_{s=0}(x)$. That is, 
\begin{equation}\label{eq:adm_tangent}
    - w\partial_x g_0 +h = \partial_s g_s|_{s=0}(x)
\end{equation} 
and
\begin{align*}
\begin{split}
\varepsilon w_{xxxx} -\lambda w_{xx} + \kappa w  & = -h \partial_x g_0  \quad \text{on } (0,1), \\
w & = 0  \quad \text{at } x = 0,1, \\
w_{xx} & = 0  \quad \text{at } x = 0,1.
\end{split}
\end{align*}
With the special inner product for the tangent space of the HV geometry~\cite[Sec.~2.2]{han2023hv}, we can express the directional derivative evaluated at the Eulerian perturbation direction $\partial_s g_s|_{s=0}(x) $ as 
\begin{align*}
   & \int_0^1 \frac{\partial J}{\partial f_0}(x) \,\partial_s g_s|_{s=0}(x) \dd x\\
    =& \lim_{s\to 0}\frac{d_{\text{HV}}^2(g_s,f_1)-d_{\text{HV}}^2(g_0,f_1)}{s}\\
   =&
   -\int_0^1 \Big(\kappa v(x,0)w(x,0)+\lambda v_x(x,0)w_x(x,0) \\
   &\quad +\varepsilon v(x,0)_{xx}w(x,0)_{xx}+z(x,0)h(x,0) \Big) \dd x.
\end{align*} 
Prop.~14 in \cite{han2023hv} implies that if $(\bar w,\bar h)$ is another pair of functions satisfying~\eqref{eq:adm_tangent}, then
\begin{align*}
  & \int_0^1 \frac{\partial J}{\partial f_0} (x)\partial_s g_s|_{s=0}(x)  \dd x \\
  %
  =& -\int_0^1\Big( \kappa v(x,0)\bar w(x,0)+\lambda v_x(x,0)\bar w_x(x,0) \\
   &\quad +\varepsilon v(x,0)_{xx}\bar w(x,0)_{xx}+z(x,0)\bar h(x,0) \Big)\dd x.
\end{align*} 
In particular, we get the following proposition by choosing $\bar w\equiv0$ and $\bar h=\partial_s g_s|_{s=0}(x) $. 
\begin{proposition}
For fixed functions $f_0,f_1\in H^1(0,1)$ let $g(s,x)=g_s(x)$ satisfy $g_0=f_0.$ If the tuple $(f,v,z)$  denotes the minimizing geodesic between $f_0$ and $f_1,$ then the derivative of $d_{\text{HV}}^2(g_s,f_1)$ with respect to $s$ at $s = 0$ can be expressed as 
$$
\lim_{s\to 0}\frac{d_{\text{HV}}^2(g_s,f_1)-d_{\text{HV}}^2(g_0,f_1)}{s}=-\int_0^1 z(x,0)\partial_s g_s|_{s=0}(x) \dd x.
$$
Based on the arbitrariness of $g_s$, we further obtain the  Fr\'echet derivative of the functional $J: f_0\rightarrow d_{\text{HV}}^2(f_0, f_1)$ given by 
\begin{equation}
  \frac{\partial J}{\partial f_0}  =   \frac{\partial d_{\text{HV}}^2(f_0,f_1)}{\partial f_0} = - z(x,0)\,.
\end{equation}
\end{proposition}

\subsection{Complex-valued HV metric and Fr\'echet derivative} Since wave data in the frequency domain is complex-valued, we extend the HV metric to complex-valued functions by matching the real and imaginary components. More precisely, if $f_0,f_1$ are complex-valued 1D functions, then for any hyperparameters $\kappa>0,\lambda\geq0$ and $\varepsilon>0$ we can define a new metric $d_{\text{HVC}}(f_0,f_1)$ as 
\begin{align*}
d_{\text{HVC}}(f_0,f_1)=\sqrt{d^2_{HV}(\Re (f_0),\Re(f_1))+d^2_{HV}(\Im(f_0),\Im(f_1))},
\end{align*} 
where $\Re(f_{j})$ and $\Im(f_j)$ denote the real and imaginary parts of the signal $f_j$, respectively.

This definition is symmetric, non-negative, and returns zero if and only if $f_0=f_1$, while the triangle inequality follows from the triangle inequality from the real-value HV metric. As in the case of real vector-valued functions, the Fr\'echet derivative of the squared HVC metric with respect to $f_0$ with $f_1$ fixed is the sum of the real and imaginary parts' Fr\'echet derivatives. This can be summarized in the following proposition.

\begin{proposition}
Let hyperparameters $\kappa>0,\lambda\geq0$ and $\varepsilon>0.$ If $f_0$ and $f_1$ are complex-valued functions $f_0,f_1:\mathbb{R}\to \mathbb{C}$, then $d_\text{HVC}(f_0,f_1)$ defines a metric. Additionaly, if $g(s,x)=g_s(x)$ satisfies $g_0=f_0$, then the derivative of $d_{\text{HVC}}^2(g_s,f_1)$ with respect to $s$ at $s = 0$ can be expressed as 
$$
\int_0^1\left( \Re(z(x,0))\Re(\theta) +i \Im(z(x,0))\Im(\theta )\right) \dd x,\,\,\theta = \partial_s g_s|_{s=0},
$$ 
where $\Re(\cdot)$ and $\Im(\cdot)$ represent a complex value's real and imaginary components, respectively. Here, $z$ is the last component of the optimal tuple connecting $f_0$ and $f_1$ under the HV geometry.
\end{proposition}

 \subsection{Parameter selection}
In \cite[Sec.~5.1]{han2023hv} the authors provided their recommendations for choosing the hyperparameters based on the scaling properties of the metric. Fig.~\ref{fig:convexpara} suggests that different values of $\kappa,\lambda$, and $\varepsilon$ highlight different properties of the metric. The exact choice of them will depend on the nature of the problem and the data provided. Generally, if the noise is significant or if our problem relies on lower-frequency data components, we want our metric to have a wider basin of attraction, which is achieved when $\varepsilon \ll1$ and $\kappa,\lambda\ll\varepsilon.$ On the other hand, if the problem requires precision and high resolution or if the initial data is close to the ground truth, we recommend choosing $\kappa,\lambda,\varepsilon\approx 1.$ In Sec.~\ref{sec:ne} for experiments $A$ and $B$, we wanted our metric to have transport-like properties, while at the same time starting from an initial velocity noticeably different from the ground truth. The hyperparameters that worked best for us are based on numerical experiments with different signals where transport-like property is required and were of the order $\kappa=10^{-10},\lambda=10^{-10},\varepsilon=10^{-7}$.

\medskip

\section{Experiments}\label{sec:ne}
\subsection{Noise-free data inversion for the full Marmousi model}\label{sec:noise_free_seismic}
The first experiment involves inverting the complete Marmousi model~\cite{martin2006marmousi2} using three objective functions: the conventional $L^2$ misfit, the trace‐by‐trace $W_2$ misfit, and the trace‐by‐trace HV misfit. By ``trace-by-trace'', we mean that the objective functions based on the $W_2$ metric and the HV metric are
\begin{align*}
J_{W_2}(c) &= \sum_{i=1}^{n_r} W_2^2(d^{\text{obs}}(x_i, t), u_k(x_i, t;c))\,,\\
J_{\text{HV}}(c) &= \sum_{i=1}^{n_r} d_{\text{HV}}^2(d^{\text{obs}}(x_i, t), u_k(x_i, t;c))\,,
\end{align*}
where $n_r$ is the total number of receivers, $c$ denotes the velocity, $d^{\text{obs}}$ is the observed wavefield while $u_k$ is the synthetic wavefield produced with the wave velocity $c$. In Appendix~\ref{sec:appendix}, we review the mathematical formulation of the $W_2$ metric, discuss the normalization of wave data required for its application, and derive the Fr\'echet derivative of the metric with respect to one of its inputs.

\begin{figure*}
\centering
\includegraphics[width=\linewidth]{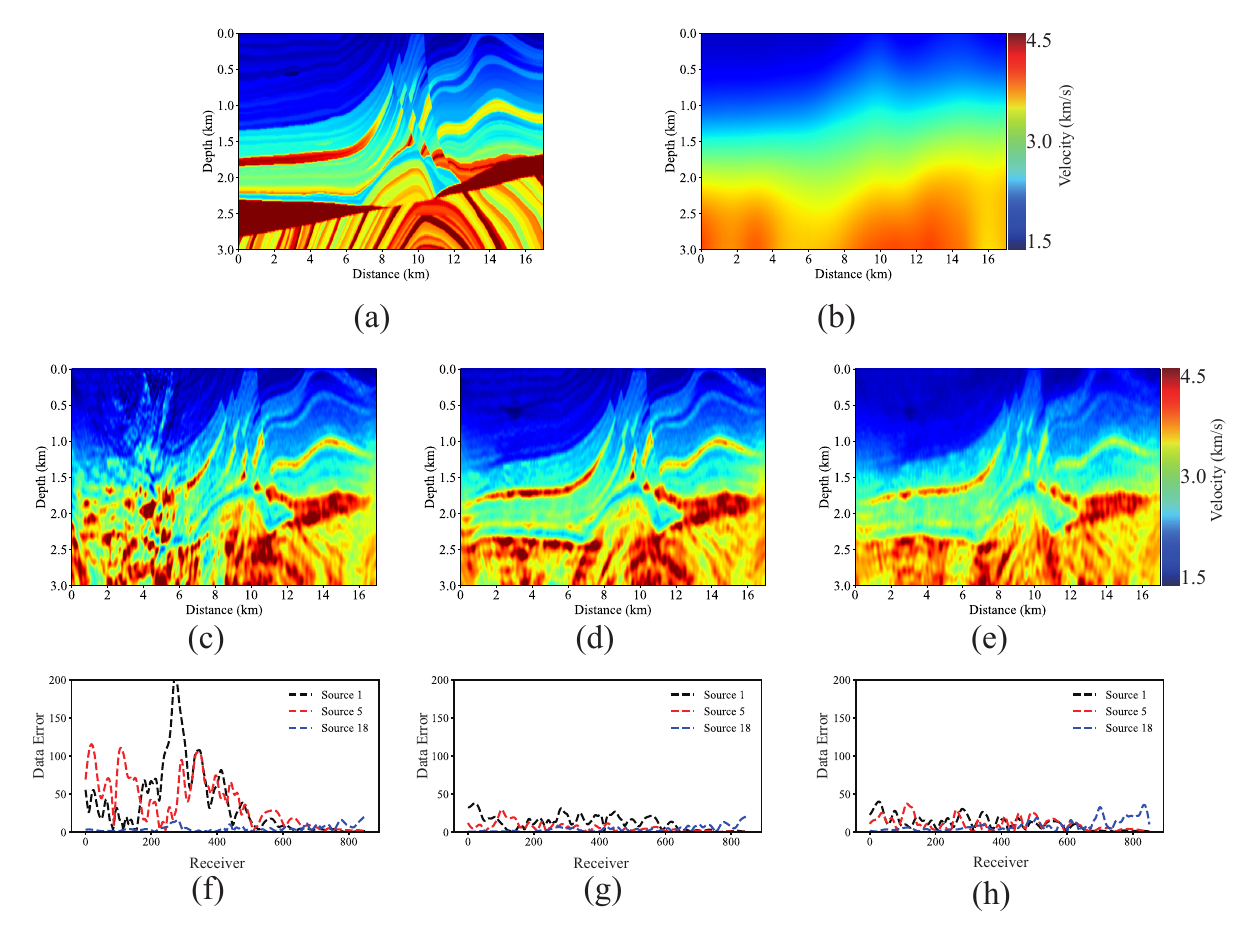}%
\caption{(a) True velocity of the Marmousi model; (b) The initial velocity in all inversions; (c)-(e) Noise-free inversion results using the $L^2$ norm, the HV metric and the $W_2$ metric, respectively; (f)-(h) The 1D data residual of sources 1, 5, 18  (the reference data subtracted by the simulated traces from velocity models in (c)-(e)) for the $L^2$ norm, the HV metric and the $W_2$ metric, respectively.}
\label{fig:mar_inf}
\end{figure*}

Fig.~\ref{fig:mar_inf}(a) shows the P-wave velocity model of the Marmousi benchmark, which spans a depth of 3~km and a width of 17~km. Twenty evenly spaced sources are positioned at a depth of 100~m in the water layer, and 851 receivers are arranged similarly to ensure comprehensive wavefield acquisition. The forward wave equation is discretized on a 20~m grid in both the $x$ and $z$ directions, with a Ricker wavelet of a peak frequency of 6~Hz serves as the source profile. The numerical solver is implemented using the 9-point optimal discretization~\cite{chen2013optimal} provided by an open-source tool~\cite{silva2019unified}, which remedies the limitations associated with the discretization and utilizes sparse LU decomposition for system inversion. Data in the 3–10~Hz frequency range is used for the inversion, and a frequency-marching strategy is adopted for sequential frequency optimization. The inversion is initiated from a blurred model obtained by smoothing the true velocity with a Gaussian filter ($\sigma = 30$; see Fig.~\ref{fig:mar_inf}(b)) and performed using the L-BFGS optimization algorithm. 

\begin{figure}[!ht]
\centering
\includegraphics[width=\linewidth]{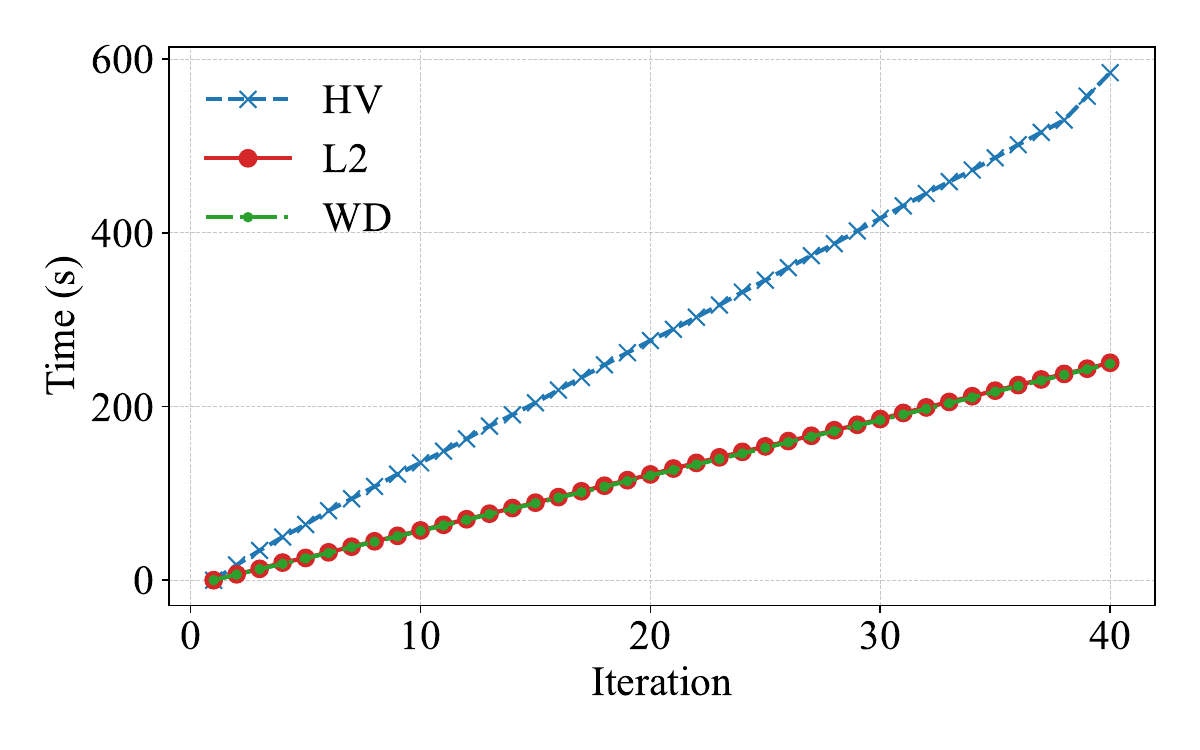}%
\caption{Computational time of 3Hz seismic inversion of the Marmousi model using three different metrics as objectives.}
\label{fig:time}
\end{figure}

Figs.~\ref{fig:mar_inf}(c)–(e) present the inversion outcomes obtained via the $L^2$ norm, the HV metric and the $W_2$ metric as the objective function and noise‐free reference wave data. The $L^2$ result exhibits numerous blocky artifacts on the left.  In contrast, the $W_2$ metric reconstructs the overall structural features of the model but reveals fragmentation and discontinuities in its high-frequency components. Meanwhile, the HV metric accurately recovers most details of the true model, yielding a smooth reconstruction that closely aligns with the ground truth. Figs.~\ref{fig:mar_inf}(f)–(h) show the data residuals among the inversion results of three metrics and the true model, further demonstrating that the HV inversion aligns most closely with the ground truth while the $L^2$ result remains trapped in a local minimum as seen in Fig.~\ref{fig:losscurve}. In Fig.~\ref{fig:time}, we present the CPU time of inverting the noise-free Marmousi model using all the metrics. We notice that on average, one iteration of HV requires roughly twice the computation time as either $L^2$ or $W_2$ due to the iterative algorithm to compute $d_\text{HV}$. 

\begin{figure*}[htbp!]
  \centering
  \includegraphics[width=\linewidth]{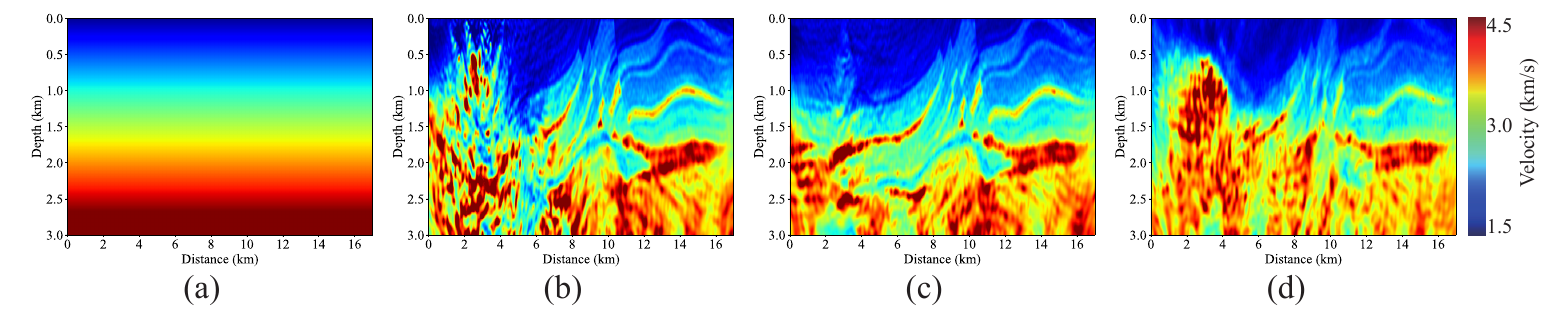}%
  \caption{(a) The initial velocity in all inversions; (b)-(d) Noise-free inversion results using the $L^2$ norm, the HV metric and the $W_2$ metric, respectively.}
  \label{fig:mar_linear}
  \end{figure*}

\begin{figure}[!ht]
  \centering
  \includegraphics[width=\linewidth]{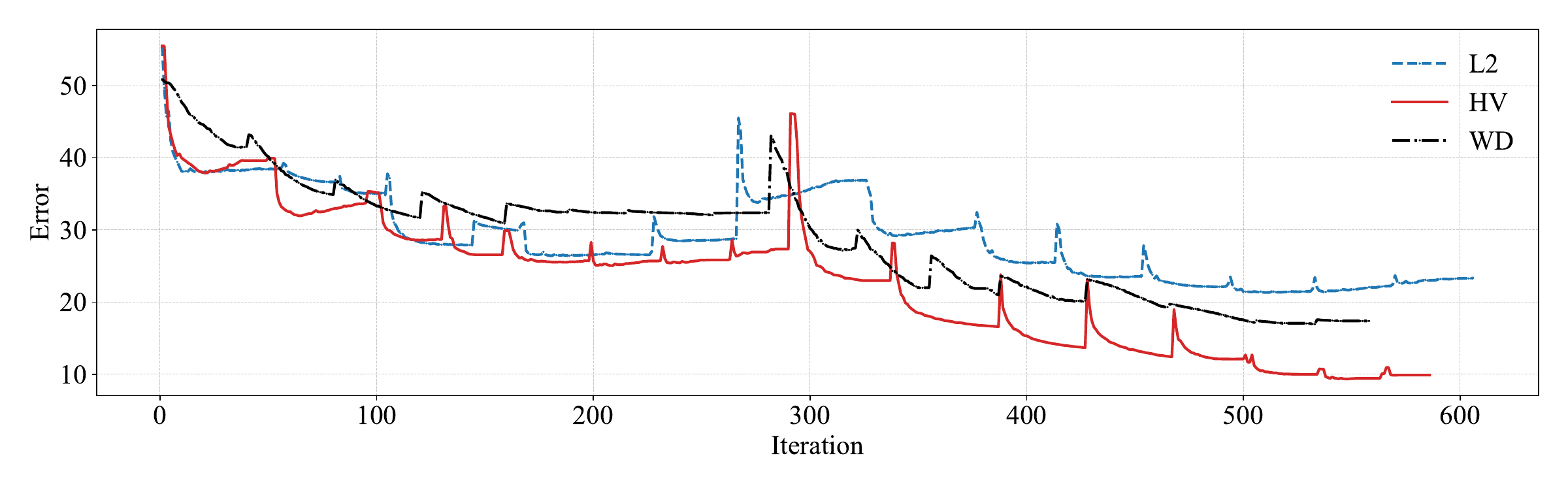}%
  \caption{Residual curves of the reconstruction results for the $L^2$, HV and $W_2$ metrics. The reconstruction employs a two-round, multi-scale frequency-marching strategy: inversion problems are solved sequentially at progressively higher wave frequencies, and the output from each round is Gaussian-filtered to serve as the initial condition for the next.}
  \label{fig:losscurve}
  \end{figure}

Fig.~\ref{fig:mar_linear} presents the inversion of the velocity model from a 1D initial model (see Fig. \ref{fig:mar_linear}(a)), representing a more challenging scenario. The HV metric successfully recovers the sound-speed information at most locations, displaying only minor artifacts on the left side. In contrast, both the $L^2$ norm and the Wasserstein distance produce extensive artifacts on the left, yielding noticeably poorer reconstruction quality.

\subsection{Noisy data inversion for the full Marmousi model}
To evaluate the robustness of the HV metric to data measurement noise, we repeat the previous experiment under more realistic conditions by adding white noise to the observed data, with a signal-to-noise (SNR) ratio of $10$~dB. 

\begin{figure*}[htbp!]
  \centering
  \includegraphics[width=\linewidth]{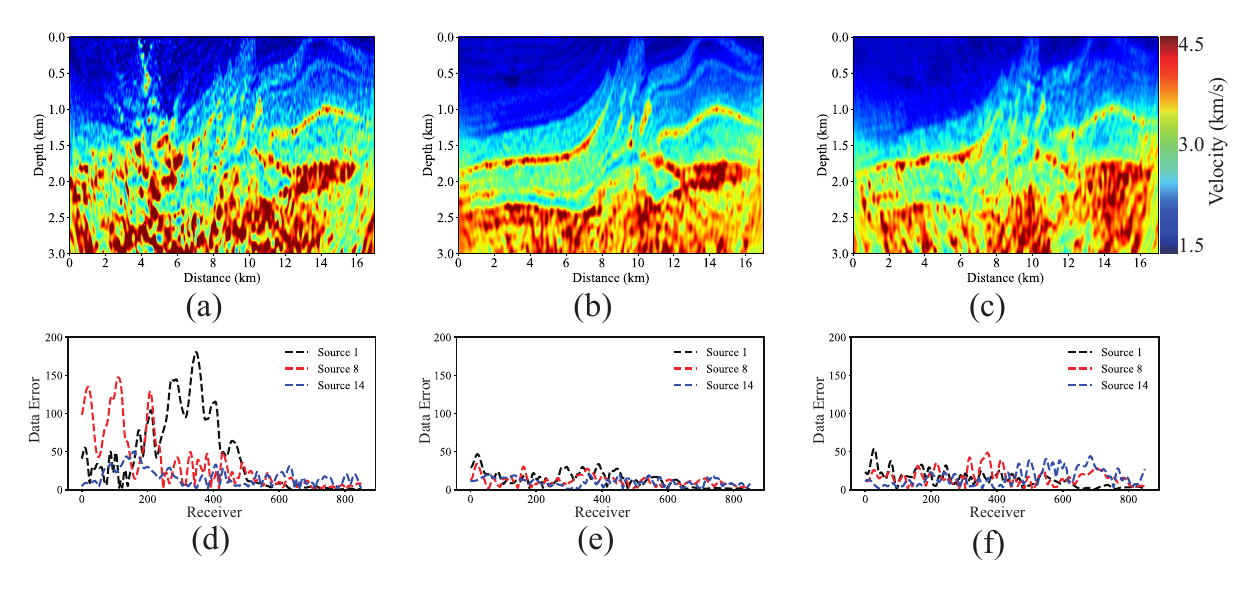}%
  \caption{Noisy data inversion result for the Marmousi benchmark (SNR is $10$~dB): (a)-(c) the inversion results using the $L^2$, $d_{\text{HV}}$ and the $W_2$ metrics, respectively; (d)-(f) The 1D data residual of sources 1, 8, 14  (the reference data subtracted by the simulated traces from velocity models in (a)-(c)) for the $L^2$ norm, the HV metric and the $W_2$ metric, respectively.}
  \label{fig:mar_10db}
  \end{figure*}
  
Fig.~\ref{fig:mar_10db}(a)–(c) present inversion results using the noisy observational data under the $L^2$, the HV, and the $2$-Wasserstein metrics, respectively. The $L^2$ metric continues to exhibit artifacts on the left side of the domain, with amplified artifacts in the upper left part, and suffers from false features. In contrast, the reconstruction obtained using the $W_2$ metric shows a marked decline in quality relative to the noise‐free case, also exhibiting  visible artifacts and discontinuities. The HV metric, however, successfully reconstructs the true model, albeit with a slight degradation in high-frequency details and continuity. The 1D data residuals shown in Fig.~\ref{fig:mar_10db}(d)–(f) corresponding to these three metrics further demonstrate that the HV metric fits the data better than the other two methods.

\subsection{Synthetic breast phantom imaging}\label{sec:simu_breast}
In this example, we demonstrate the advantages of the HV metric in USCT breast imaging. The virtual annular USCT system consisted of a circular transducer array of 220~mm diameter ring along which 256 transducers were equispaced and acted as receivers. During data acquisition, one transducer is selected to emit an excitation pulse while the remaining sensors record the responses; this process is repeated $256$ times. We obtained $25$ synthetic breast wave speed models, including all four common breast tissue types (all‐fatty, fibroglandular, heterogeneous, and extremely dense), from an open‐source 2D breast wave speed model dataset~\cite{zeng2025openwaves}. In our simulations, the Convergent Born Series (CBS)~\cite{osnabrugge2016convergent,aghamiry2022large} was used to simulate both the observed data (300~kHz frequency data) and FWI imaging. CBS incorporates a built-in preconditioner into the Born expansion,
guaranteeing convergence even in strongly scattering, heterogeneous media. The initial velocity model is defined as a water‐filled medium, indicating that no prior information is available for reconstruction. To simulate realistic measurement conditions, additive Gaussian noise was introduced, which produced an SNR of 10~dB.

\begin{figure*}[htbp!]
\centering
\includegraphics[width=0.8\linewidth]{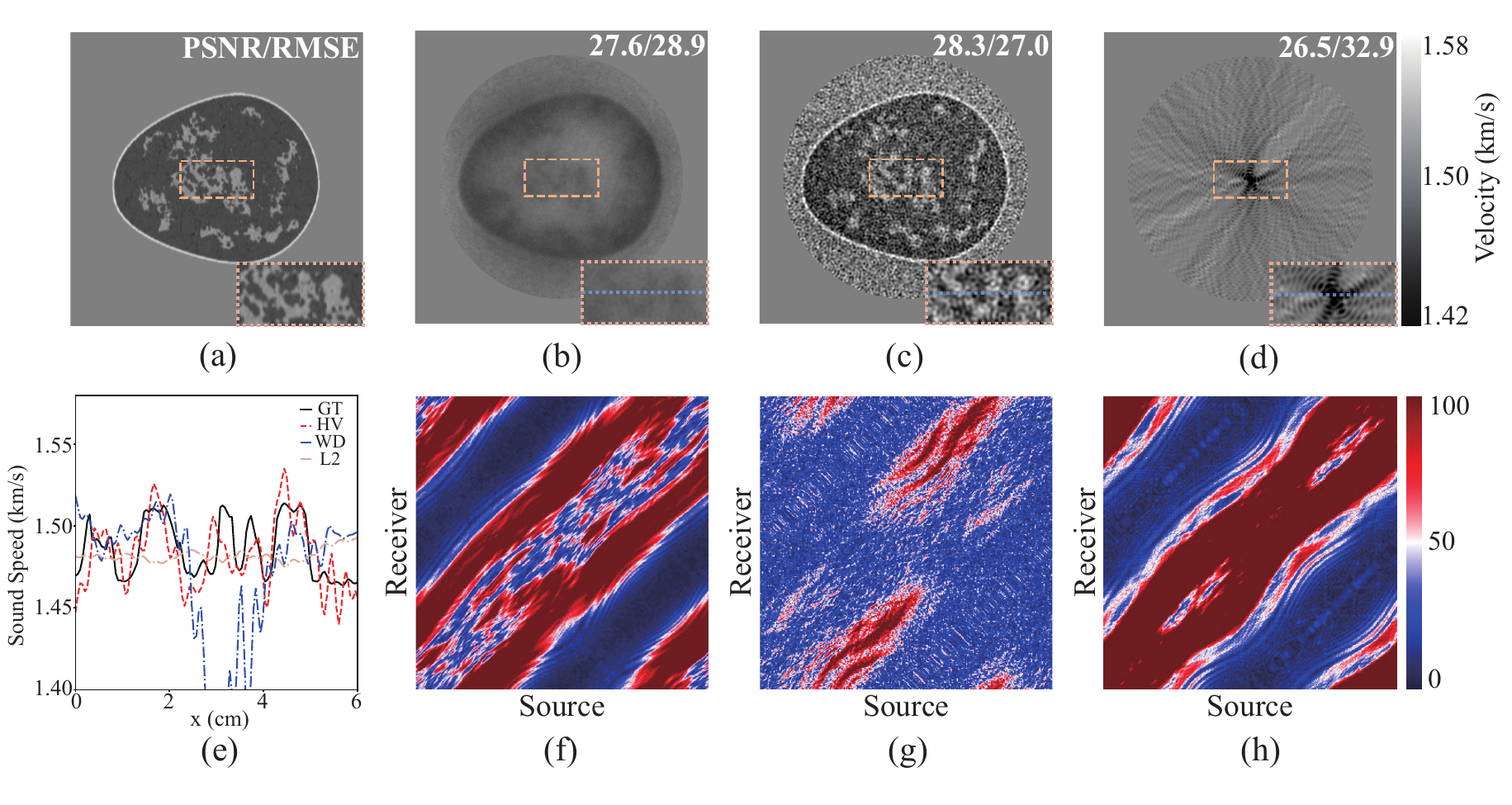}%
\caption{Synthetic breast phantom imaging: (a)  True synthetic breast velocity model; (b)-(d) The inversion results using the $L^2$, $d_{\text{HV}}$ and the $W_2$ metrics, respectively;  (e)~1D line profile of the reconstruction velocity model using $L^2$, $d_{\text{HV}}$, and $W_2$ metrics, respectively; (f)-(h) the relative data residual at final inversion results using the $L^2$, $d_{\text{HV}}$ and the $W_2$ metrics, respectively. }
\label{fig:syn_breast}
\end{figure*}

The reference wave speed model is presented in Fig.~\ref{fig:syn_breast}(a). Figs.~\ref{fig:syn_breast}(b–d) display the reconstruction results of a heterogeneous breast obtained using three misfit functions, $L^2$, $d_\text{HV}$, and the $W_2$ metric, while Fig.~\ref{fig:syn_breast}(e) provides one‑dimensional velocity profiles at selected locations to offer a more detailed comparison. The inversion results indicate that the HV metric-based reconstructions accurately captured detailed breast structures, including skin, fat, and glandular tissue, despite the presence of noise-induced imaging artifacts. In contrast, $L^2$ metric-based FWI is trapped in local minima and fails to recover the true wave speed distribution, whereas the $W_2$ metric does not provide an effective gradient and consequently fails to yield a meaningful reconstruction. The line profile in Fig.~\ref{fig:syn_breast}(e) and the data residual in Figs.~\ref{fig:syn_breast}(f-h) further validate the degree of alignment between the inversion result and the true velocity model.

\begin{figure*}[ht!]
\centering
\includegraphics[width=0.8\linewidth]{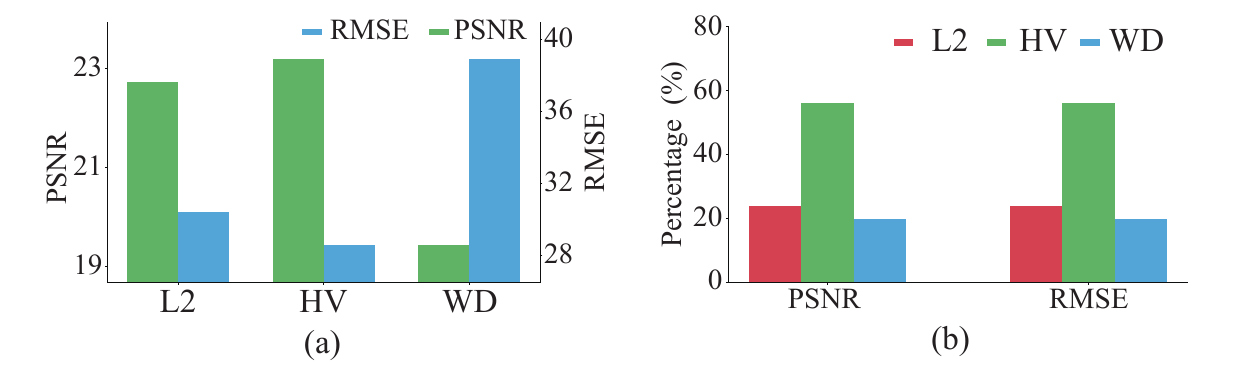}%
\caption{Synthetic breast phantom imaging: (a) Average PSNR and RMSE of reconstruction result using $L^2$, HV, and $W_2$ metrics as the objective functions; (b) The percentage of experiments in which a model outperformed the others.}
\label{fig:metric_syn_breast}
\end{figure*}

To quantify the imaging quality of the breast phantom, we evaluated the root mean square error (RMSE) and peak signal-to-noise ratio (PSNR) between 25 breast reconstruction results and their respective reference models, as shown in Fig.~\ref{fig:metric_syn_breast}. The comparison demonstrates that reconstructions based on the HV metric outperform the other two metrics in terms of both RMSE and PSNR and yield the best quality in over half of the samples.

\subsection{Open-source in vivo breast USCT model}
\begin{figure*}[htbp!]
\centering
\includegraphics[width=0.8\linewidth]{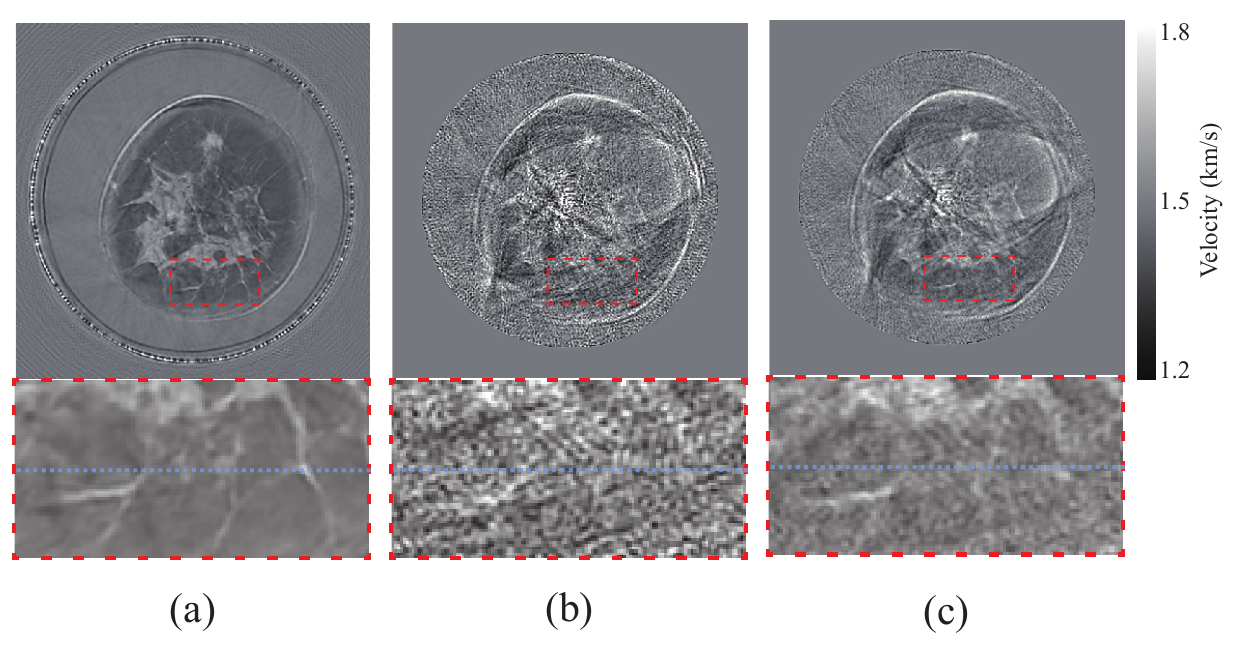}%
\caption{Open-source in-vivo breast USCT model: (a) Reference velocity model; (b) Reconstruction result using the $L^2$ metric; (c) Reconstruction result using the HV metric.}
\label{fig:exp_breast}
\end{figure*}

To experimentally validate the HV metric, we test it and the $L^2$ norm on an open-source \textit{in vivo} breast USCT dataset from the Karmanos Cancer Institute (IRB Approval No.~040912M1)~\cite{ali2024ringFWI2D} where we utilized the dataset corresponding to breasts with malignant lesions. To assess the robustness of the HV metric, we excluded low‑frequency components and restricted the frequency span for wave‑speed reconstruction to [0.40, 0.45, 0.50, 0.55, 0.60, 0.65, 0.70] MHz, thereby creating a more challenging imaging scenario. The numerical solver and initial velocity model were identical to those described in Sec.~\ref{sec:simu_breast}. Since we do not know the ground truth for this real-world application, the reference velocity model is given by a reconstruction result with $20$ frequencies data~\cite{ali2024ringFWI2D}, as shown in Fig.~\ref{fig:exp_breast}. 

As demonstrated in Fig.~\ref{fig:exp_breast}, the $L^2$ metric yields reconstructions contaminated by pronounced noise-induced artifacts, which obscure critical anatomical features such as glandular tissues within the breast. In contrast, the HV metric successfully resolves structural details despite utilizing only partial frequency-band data, with significantly reduced noise-related blurring compared to the $L^2$-based results. This comparative analysis demonstrates the potential of the HV metric in real-world imaging scenarios with imperfect observational data.

\section{Conclusion}\label{sec:conclusion}
In this study, we introduce a novel misfit function, the HV metric ($d_{\text{HV}}$), for frequency-domain full-waveform inversion (FWI), specifically designed to address the cycle-skipping problem and improve robustness to noisy measurements. The HV metric accommodates discontinuities in the underlying signals and jointly captures horizontal and vertical deformations, inheriting key advantages of the Wasserstein metric. By extending the original HV metric into a trace-by-trace formulation for complex-valued signals, we apply it successfully across a range of FWI tasks. The HV metric consistently outperforms both the conventional $L^2$ norm and the quadratic Wasserstein metric ($W_2$) in seismic inversion and ultrasound computed tomography (USCT), demonstrating superior robustness to observational noise. Additionally, inversion results from real clinical USCT data highlight the effectiveness of the HV metric even at high frequencies.

Our results indicate that $d_{\text{HV}}$ is a highly promising candidate for a misfit function in frequency‑domain FWI. Despite the advantages outlined above, several limitations warrant further investigation. First, although our numerical examples reveal the noise robustness of $d_{\text{HV}}$, this phenomenon lacks theoretical justification. Moreover, since the computation of $d_{\text{HV}}$ requires solving a PDE‑constrained optimization problem, its computational cost increases markedly with the number of transducers, which poses a major challenge for its application to large-scale three-dimensional imaging tasks. In future work, we plan to investigate the noise robustness behavior of $d_{\text{HV}}$ via the perturbation analysis. By refining the numerical schemes, initialization strategies, and optimization algorithms tailored for the $d_{\text{HV}}$-based inversion, we expect not only to improve its efficiency but also ultimately enable its successful application to real‑world 3D FWI-related problems.  

\section*{Acknowledgments}
M.N.~and Y.Y.~were supported by the National Science Foundation through grant DMS-2409855 and the Office of Naval Research through grant N00014-24-1-2088. This work was done in part while Z.Z.~was visiting Cornell University in Fall 2024. We are very grateful to Prof. He Sun and Prof Yubing Li for sharing the \textit{OpenWaves} Synthetic Breast Dataset and the USCT algorithm and to Prof. Dejan Slep\v{c}ev for providing the Fr\'echet derivative formula of the HV metric.

\section{Appendix}
\subsection{The quadratic Wasserstein metric}\label{sec:appendix}
In this appendix, we review the mathematical formulation of the $2$-Wasserstein metric ($W_2$), discuss the normalization of wave data required for its application, and derive the Fr\'echet derivative of the metric with respect to one of its inputs. For simplicity, we assume the input data are one-dimensional signals.

\subsubsection{Formulation}
Mathematically, let $f_0:\mathbb{R}\rightarrow\mathbb{R}^{+}$ and $f_1:\mathbb{R}\rightarrow\mathbb{R}^{+}$ be two probability density functions on $\mathbb{R}$ that satisfy $\int_{\mathbb{R}}f_0(s)ds = \int_{\mathbb{R}}f_1(s)ds = 1$, then the mathematical definition of the quadratic Wasserstein
metric between $f_0$ and $f_1$ will be
\begin{equation}
    W_2(f_0, f_1) = \sqrt{\inf_{T \in \mathcal{M}} \int_{\mathbb{R}} |s - T(s)|^2 f_0(s) \, ds},
\end{equation}
in which $\mathcal{M}$ is the set of all the rearrange maps from $f_0$ to $f_1$. According to the optimal transportation theorem for a quadratic cost on $\mathbb{R}, $~\cite{villani2021topics}, the optimal transportation cost and the optimal map are given by
\begin{align}
  W_2(f_0, f_1) &= \sqrt{\int_0^1 \left| F^{-1}(s) - G^{-1}(s) \right|^2 \, ds},\\
  T(s) &= G^{-1}(F(s)),
\end{align}
where $F(s)=\int_{-\infty}^{s} f_0(\tau) \, d\tau$ and $G(s)= \int_{-\infty}^{s} f_1(\tau) \, d\tau$ are the corresponding cumulative distribution functions of the probability density functions $f_0(s)$ and $f_1(s)$, respectively.

\subsubsection{Data normalization and the Fr\'echet derivative}
In frequency-domain FWI, the complex-valued measurements cannot be directly regarded as distributions, presenting challenges for computing the Wasserstein metric. Unlike in the time domain, where the signal at individual receivers is treated as a one-dimensional time series, in the frequency domain, it is more appropriate to consider the data obtained from all receivers in response to a single excitation as the boundary measurement.

To satisfy the requirements of the Wasserstein metric for mass conservation and nonnegative, scalar-valued data, we treat the complex-valued seismic measurements from the frequency domain as two independent one-dimensional real-valued time series. The two signals are then normalized by the operator $\mathcal{T}: \mathbb{R}\rightarrow \mathbb{R}^{+}$, resulting in the modified misfit function
\[
\begin{aligned}
d_{W_2}(f_0,f_1) &=& W_2^2\left(\frac{\mathcal{T}(\Re{f_0})}{\langle \mathcal{T}(\Re{f_0} \rangle },\frac{\mathcal{T}(\Re{f_1})}{\langle \mathcal{T}(\Re{f_1}) \rangle}\right) \\
&+& W_2^2\left(\frac{\mathcal{T}(\Im{f_0})}{\langle \mathcal{T}(\Im{f_0})  \rangle},\frac{\mathcal{T}(\Im{f_1})}{\langle \mathcal{T}(\Im{f_1}) \rangle}\right)
\end{aligned}
\] 
where operators $\Re$ and $\Im$ extract the real and imaginary parts of the measured data, and the operator $\langle \cdot \rangle$ denotes the integral over the real axis, i.e.,
\[
\langle f \rangle = \int_{a}^{b} f(s)ds
\]
where we assume that the signal is supported over $[a,b]$ for simplicity.

To obtain the gradients for the current velocity model, the adjoint-state method is the most commonly employed approach. A critical component of this method is the derivation of the Fr\'echet derivative of the objective function, which plays the role of the source term for the adjoint equation. Let $\delta f_0$ be a small perturbation of $f_0$. Then the perturbation to the squared Wasserstein metric ($W_2^2(f_0,f_1)$), denoted by $\delta \mathcal{W}$, is given by
\[
  \delta \mathcal{W} = \int_0^{t_f} \left(\int_{0}^{s} 2(\tau - T(\tau))d\tau\right)  \delta f_0 (s) \, ds.
\]
Next, for a differentiable and one-to-one normalization operator $\mathcal{T}$, the Fr\'echet derivative of the resulting signal $\mathcal{S}(f_0) = \frac{\mathcal{T}(f_0)}{\langle \mathcal{T}(f_0)\rangle}$ with respect to $f_0$ is then given by
\[
\delta\mathcal{S} = \frac{\mathcal{T}'(f_0)\delta f_0}{\langle \mathcal{T}(f_0)\rangle} - \frac{\mathcal{T}(f_0)\langle \mathcal{T}'(f_0)\delta f_0\rangle}{\langle \mathcal{T}(f_0)\rangle^2}\,.
\]

To sum up, the Fr\'echet derivative  of $\nabla_f d_{W_2}(f_0,f_1)$ with respect to the original signal $f_0$ with $f_1$ fixed can be obtained by the chain rule
\[
\begin{aligned}
& \nabla_{f_0} d_{W_2}(f_0,f_1) \\
=&  \ \nabla_{\mathcal{S}(\text{R}(f_0))} W_2^2(\mathcal{S}(\text{R}(f_0)),\mathcal{S}(\text{R}(f_1))) \nabla_{\text{R}(f_0)} \mathcal{S}(\text{R}f_0) \\
 & + \ i \nabla_{\mathcal{S}(\text{I}(f_0))} W_2^2(\mathcal{S}(\text{I}(f_0)),\mathcal{S}(\text{I}(f_1))) \nabla_{\text{I}(f_0)} \mathcal{S}(\text{I}f_0).
\end{aligned}
\]
In this work, we use the linear transformation to define our operator $\mathcal{T}$ following~\cite{yang2018application}.

\bibliographystyle{IEEEtran}
\bibliography{ref.bib}

\end{document}